\documentclass{amsproc}
\usepackage{graphicx}
\usepackage{amscd}
\usepackage{amsmath}
\usepackage{amsfonts}
\usepackage{amssymb}
\newtheorem{theorem}{Theorem}
\theoremstyle{plain}
\newtheorem{acknowledgement}{Acknowledgement}

\newtheorem{definition}{Definition}
\newtheorem{example}{Example}

\newtheorem{proposition}{Proposition}
\newtheorem{remark}{Remark}

\numberwithin{equation}{section}

\begin{document}
\title[A solution of one J. Bourgain's problem]{A solution of one J. Bourgain's problem}
\author{Eugene Tokarev}
\address{B.E. Ukrecolan, 33-81 Iskrinskaya str., 61005, Kharkiv-5, Ukraine}
\email{tokarev@univer.kharkov.ua}
\dedicatory{Dedicated to the memory of S. Banach.}
\begin{abstract}It is proved that there exists a separable reflexive Banach space $W$ that
contains an isomorphic image of every separable superreflexive Banach space.
This gives the answer on one J. Bourgain's question
\end{abstract}
\maketitle

\section{Introduction}

In 1980 J. Bourgain [1] posed a question:

\textit{Whether there exists a separable reflexive Banach space universal for
all separable superreflexive Banach spaces?}

Recall that a Banach space $X$ is said to be \textit{universal} for a class
$\mathcal{K}$ of Banach spaces if every member of $\mathcal{K}$ is isomorphic
to a subspace of $X$.

A partial solution of the problem was obtained in [2], where the existence of
a separable reflexive Banach space that contains all separable superreflexive
Banach spaces, which enjoy the approximation property, was proved.

In the article the existence of the corresponding universal space will be
shown without additional restrictions.

The idea of the construction is rather simple.

Let $\mathcal{B}$ be a (proper) class of all Banach spaces; $1<p\leq2\leq
q<\infty$. Let $C(p,q,N)$ $\subset\mathcal{B}$ be a set of all separable
superreflexive Banach spaces that are both of type $p$ and of cotype $q$, with
corresponding constants of type and cotype that do not exceed $N$.

It may be organized a Banach space $X$ (certainly, non-separable) as an
$l_{2}$-sum of all spaces from $C(p,q,N)$:%
\[
X=(\sum\oplus\{E:E\in C(p,q,N)\})_{2},
\]
where $\left(  \sum\oplus\{X_{i}:i\in I\}\right)  _{2}$ is a Banach space of
all families $\{x_{i}\in X_{i}:i\in I\}=\frak{x}$ with the finite norm%
\[
\left\|  \frak{x}\right\|  _{2}=\sup\{(\sum\{\left\|  x_{i}\right\|  _{X_{i}%
}^{2}:i\in I_{0}\})^{1/2}:I_{0}\subset I;\text{ }\operatorname{card}\left(
I_{0}\right)  <\infty\}.
\]

This space generates a class $X^{f}$ of all Banach spaces that are
\textit{finitely equivalent }to $X$, i.e. of those Banach spaces $Y$, which
are finitely representable in $X$ (shortly, $Y<_{f}X$) and in which the space
$X$ is finitely representable as well. Thus,%
\[
X^{f}=\{Y\in\mathcal{B}:Y<_{f}X\text{ \ and \ }X<_{f}Y\}.
\]

Recall (cf. [3]) that a Banach space $E$ is said to be an \textit{isomorphic
envelope }of a class $E^{f}$ provided any Banach space $Z$, which is finitely
representable in $E$ and whose dimension $\dim(Z)$ does not exceed $\dim(E)$
is isomorphic to a subspace of $E$. Certainly, there no guarantee that $X^{f}$
contains a space, which is the separable isomorphic envelope of this class.

The main step of the construction is to built a new class $\left(
W(X)\right)  ^{f}$ (certainly, depending on $X$) such that

\begin{itemize}
\item  Any space of $\left(  W(X)\right)  ^{f}$ has the same type and cotype
as $X$;

\item $X$ is finitely representable in $W(X)$;

\item $W(X)$ is superreflexive provided $X$ has a such property;

\item  The class $\left(  W(X)\right)  ^{f}$ contains a separable space (say,
$E_{W}$) which is an isomorphic envelope of $\left(  W(X)\right)  ^{f}$.
\end{itemize}

Next steps of the construction are obvious. For every $n\in\mathbb{N}$
consider a class $C_{n}=C(1-n^{-1},n,n)$ (described above) and a corresponding
space
\[
X_{n}=(\sum\oplus\{E:E\in C_{n}\})_{2}.
\]
Let $\left(  W_{n}\right)  ^{f}=\left(  W(X_{n})\right)  ^{f}$ be the
corresponding class. Let $E_{n}$ be a separable isomorphic envelope of this
class. Clearly,
\[
G=(\sum\nolimits_{n<\infty}\oplus E_{n})_{2}%
\]
is a desired space. Indeed, $G$ is reflexive because of it is the $l_{2}$-sum
of reflexive spaces. $G$ is separable since it is a countable sum of separable
spaces. At least, if $Z$ is a superreflexive separable Banach space then $Z$
is isometric to a subspace of one of $X_{n}$'s (by the definition of $X_{n}$).
Therefore $Z$ is finite representable in a corresponding space $W_{n}$ and,
hence, is isomorphic to a subspace of $E_{n}$ (by the definition of an
isomorphic envelope). Consequently, $Z$ is isomorphic to a subspace of $G$.

Below it will be shown that such a space $G$ is not unique: there exists a
continuum number of pairwise non-isomorphic spaces $G_{\alpha}$ with the same property.

\section{Definitions and notations}

\begin{definition}
Let $X$, $Y$ are Banach spaces, $\lambda<\infty$. $X$ is $\lambda
$-\textit{finitely representable} in $Y$  if for  every finite dimensional
subspace $A$ of $X$ there exists a subspace $B$ of $Y$ and an isomorphism
$u:A\rightarrow B$ such that $\left\|  u\right\|  \left\|  u^{-1}\right\|
<\lambda$.

$X$ is finite representable in $Y$\ (in symbols: $X<_{f}Y$) if $X$ is
$\lambda$-\textit{finitely representable} in $Y$ \ for every $\lambda>1$.

$X$ is crudely finite representable in $Y$ if $X$ is $\lambda$%
-\textit{finitely representable} in $Y$ \ for some $\lambda>1.$
\end{definition}

Notice \ that $X$ is crudely finite representable in $Y$ if and only if $X$ is
isomorphic to a space $X_{1}$, which is finitely representable in $Y$.

For any two Banach spaces $X$, $Y$ their \textit{Banach-Mazur distance }is
given by%
\[
d(X,Y)=\inf\{\left\|  u\right\|  \left\|  u^{-1}\right\|  :u:X\rightarrow
Y\},
\]
where $u$ runs all isomorphisms between $X$ and $Y$ and is assumed, as usual,
that $\inf\varnothing=\infty$.

It is well known that $\log d(X,Y)$ forms a metric on each class of isomorphic
Banach spaces, where almost isometric Banach spaces are identified.

Recall that Banach spaces $X$ and $Y$ are \textit{almost isometric} if
$d(X,Y)=1$. Surely, any almost isometric finite dimensional Banach spaces are isometric.

The set $\frak{M}_{n}$ of all $n$-dimensional Banach spaces, equipped with
this metric, is the compact metric space, called \textit{the Minkowski
compact} $\frak{M}_{n}$.

The disjoint union $\cup\{\frak{M}_{n}:n<\infty\}=\frak{M}$ is a separable
metric space, which is called the \textit{Minkowski space}.

Consider a Banach space $X$. Let $H\left(  X\right)  $ be a set of all its
different finite dimensional subspaces (isometric finite dimensional subspaces
of $X$ in $H\left(  X\right)  $ are identified). Thus, $H\left(  X\right)  $
may be regarded as a subset of $\frak{M}$, equipped with a restriction of the
metric topology of $\frak{M}$.

Of course, $H\left(  X\right)  $ need not to be a closed subset of $\frak{M}$.
Its closure in $\frak{M}$ will be denoted by $\overline{H\left(  X\right)  }$.
From definitions it follows that $X<_{f}Y$ if and only if $\overline{H\left(
X\right)  }\subseteq\overline{H\left(  Y\right)  }$.

\begin{definition}
Spaces $X$ and $Y$ are \textit{finitely equivalent }(in symbols: $X\sim_{f}Y$)
if $X<_{f}Y$ and $Y<_{f}X$.
\end{definition}

Therefore, $X\sim_{f}Y$ if and only if $\overline{H\left(  X\right)
}=\overline{H\left(  Y\right)  }$. Thus, there is a one to one correspondence
between classes of finite equivalence
\[
X^{f}=\{Y\in\mathcal{B}:X\sim_{f}Y\}
\]
and closed subsets of $\frak{M}$ of kind $\overline{H\left(  X\right)  }$.
Indeed, all spaces $Y$ from $X^{f}$ have the same set $\overline{H\left(
X\right)  }$. This set, uniquely determined by $X$ (or, equivalently, by
$X^{f}$), will be denoted $\frak{M}(X^{f})$ and will be referred to as
\textit{the Minkowski's base of the class} $X^{f}$.

The following classification of Banach spaces due to L. Schwartz [4].

\begin{definition}
For a Banach space $X$ its $l_{p}$-\textit{spectrum }$S(X)$ is given by%
\[
S(X)=\{p\in\lbrack0,\infty]:l_{p}<_{f}X\}.
\]
\end{definition}

Certainly, if $X\sim_{f}Y$ then $S(X)=S(Y)$. Thus, the $l_{p}$-spectrum $S(X)$
may be regarded as a property of the whole class $X^{f}$. So, notations like
$S(X^{f})$ are of obvious meaning.

\begin{definition}
Let $X$ be a Banach space. It is called:

\begin{itemize}
\item $c$-\textit{convex,} if $\infty\notin S(X)$;

\item $B$-\textit{convex,} if $1\notin S\left(  X\right)  $;

\item \textit{Finite universal,} if $\infty\in S(X)$.
\end{itemize}
\end{definition}

As it was shown in [5], the $l_{p}$-spectrum is closely connected with notions
of type and cotype. Recall the definition.

\begin{definition}
Let $1\leq p\leq2\leq q\leq\infty$. A Banach space $X$ is said to be of
\textit{type} $p$, respectively, of \textit{cotype }$q$, if for every finite
sequence $\{x_{n}:n<N\}$ of its elements
\[
\int_{0}^{1}\left\|  \sum\nolimits_{n=0}^{N-1}r_{n}\left(  t\right)
x_{n}\right\|  dt\leq t_{p}\left(  X\right)  \left(  \sum\nolimits_{n=0}%
^{N-1}\left\|  x_{n}\right\|  ^{p}\right)  ^{1/p},
\]
respectively,%
\[
\left(  \sum\nolimits_{n=0}^{N-1}\left\|  x_{n}\right\|  ^{q}\right)
^{1/q}\leq c_{q}\left(  X\right)  \int_{0}^{1}\left\|  \sum\nolimits_{n=0}%
^{N-1}r_{n}\left(  t\right)  x_{n}\right\|  dt,
\]
where $\left\{  r_{n}\left(  t\right)  :n<\infty\right\}  $ are Rademacher functions.
\end{definition}

When $q=\infty$, the sum $\left(  \sum_{n=0}^{N-1}\left\|  x_{n}\right\|
^{q}\right)  ^{1/q}$ must be replaced with $\sup\nolimits_{n<N}\left\|
x_{n}\right\|  $. Constants $t_{p}\left(  X\right)  $ and $c_{q}\left(
X\right)  $\ in these inequalities depend only on $X$. Their least values
$T_{p}\left(  X\right)  =\inf t_{p}\left(  X\right)  $ and $C_{q}\left(
X\right)  =\inf c_{q}\left(  X\right)  $, are called \textit{the type}
$p${\small -}\textit{constant} $T_{p}\left(  X\right)  $ and \textit{the
cotype} $q$\textit{-constant }$C_{q}\left(  X\right)  $.

Every Banach space is of type $1$ and of cotype $\infty$.

If $X$ is of type $p$ and of cotype $q$ with the constants $T_{p}\left(
X\right)  =T$, $C_{q}\left(  X\right)  =C$ than any $Y\in X^{f}$ is of same
type and cotype and its type-cotype constants are equal to those of $X$. Thus,
it may be spoken about the type and cotype of the whole class $X^{f}$ of
finite equivalence. Notice that $X$ and $Y$ are isomorphic then these spaces
are of the same type and cotype.

It is known (see [5]) that
\begin{align*}
\inf S(X)  &  =\sup\{p\in\lbrack1,2]:T_{p}\left(  X\right)  <\infty\};\\
\sup S(X)  &  =\inf\{q\in\lbrack2,\infty]:C_{q}\left(  X\right)  <\infty\}.
\end{align*}

\begin{definition}
A Banach space $X$ is said to be \textit{superreflexive} if every space of the
class $X^{f}$ is reflexive.
\end{definition}

Equivalently, $X$ is superreflexive if any $Y<_{f}X$ is reflexive. Clearly,
any superreflexive Banach space is $B$-convex. Since $S(X)$ is either
$[1,\infty]$ or is a closed subset of $[1,\infty)$ (cf. [6]), any
superreflexive space is of \textit{non-trivial} (i.e., non equal to $1$) type.

\begin{definition}
(Cf. $[7]$). Let $X$ be a Banach space; $Y$ - its subspace. $Y$ is said to be
a reflecting subspace of $X$ (symbolically: $Y\prec_{u}X$) if for every
$\varepsilon>0$ and every finite dimensional subspace $A\hookrightarrow X$
there exists an isomorphic embedding $u:A\rightarrow Y$ such that $\left\|
u\right\|  \left\|  u^{-1}\right\|  \leq1+\varepsilon$ and $u\mid_{A\cap
Y}=Id_{A\cap Y}$.
\end{definition}

As it was shown in [7], if $Y\prec_{u}X$ then $Y^{\ast\ast}$ is an image of a
norm one projection $P:X^{\ast\ast}\rightarrow Y^{\ast\ast}$ (under canonical
embedding of $Y^{\ast\ast}\ $into $X^{\ast\ast}$).

\begin{definition}
(Cf. $[8]$). A Banach space $E$ is said to be existentialy closed in a class
$X^{f}$ if for any isometric embedding $i:E\rightarrow Z$ into an arbitrary
space $Z\in X^{f}$ its image $iE$ is a reflecting subspace of $Z$:
$iY\prec_{u}Z$.
\end{definition}

A class of all spaces $E$ that are existentialy closed in $X^{f}$ is denoted
by $\mathcal{E}\left(  X^{f}\right)  $. In [8] it was shown that for any
Banach space $X$ the class $\mathcal{E}\left(  X^{f}\right)  $ is nonempty;
moreover, any $Y<_{f}X^{f}$ may be isometricaly embedded into some
$E\in\mathcal{E}\left(  X^{f}\right)  $ of the dimension $\dim(E)=\max
\{\dim(Y),\omega\}$($\omega$ denotes the first infinite ordinal number).

\section{Quotient closed divisible classes of finite equivalence}

In this section it will be shown how to enlarge a Minkowski's base
$\frak{M}(X^{f})$ of a certain $B$-convex (resp., superreflexive) class
$X^{f}$ to obtain a set $\frak{N}$ which will be a Minkowski's base
$\frak{M}(W^{f})$ for some class $W^{f}$, which holds the $B$-convexity (resp.
superreflexivity) and type and cotype of $X^{f}$, which, as it will be shown
in the next section, contains an isomorphic separable envelope.

\begin{definition}
A class $X^{f}$ (and its Minkowski's base $\frak{M}(X^{f})$) is said to be
divisible if some (equivalently, any) space $Z\in X^{f}$ is finitely
representable in any its subspace of finite codimension.
\end{definition}

\begin{example}
Any Banach space $X$ may be isometricaly embedded into a space
\[
l_{2}(X)=(\sum\nolimits_{i<\infty}\oplus X_{i})_{2},
\]
where all $X_{i}$'s are isometric to $X$. Immediately, $l_{2}(X)$ generates a
divisible class $\mathsf{D}_{2}(X^{f})=\left(  l_{2}(X)\right)  ^{f}$ which is
of the same type and cotype as $X^{f}$ and is superreflexive if and only if
$X^{f}$ is superreflexive.
\end{example}

To distinguish between general divisible classes and classes of type
$\mathsf{D}_{2}(X^{f})$, the last ones will be called $\mathit{2}%
$\textit{-divisible classes}.

\begin{definition}
A class $X^{f}$ (and its Minkowski's base $\frak{M}(X^{f})$) is said to be
quotient closed if for any $A\in\frak{M}(X^{f})$ and its subspace
$B\hookrightarrow A$ the quotient $A/B$ belongs to $\frak{M}(X^{f})$.
\end{definition}

Let $K\subseteq\frak{M}$ be a class of finite dimensional Banach spaces
(recall, that isometric spaces are identified). Define operations $H$, $Q$ and
$^{\ast}$ that transform a class $K$ to another class of finite dimensional
Banach spaces - $H(K)$; $Q(K)$ or $(K)^{\ast}$ respectively. Namely, let%
\begin{align*}
H(K)  &  =\{A\in\frak{M}:A\hookrightarrow B;\text{ \ }B\in K\}\\
Q(K)  &  =\{A\in\frak{M}:A=B/F;\text{ \ }F\hookrightarrow B;\text{ \ }B\in
K\}\\
(K)^{\ast}  &  =\{A^{\ast}\in\frak{M}:A\in K\}
\end{align*}

In words, $H(K)$ consists of all subspaces of spaces from $K$; $Q(K)$ contains
all quotient spaces of spaces of $K$; $(K)^{\ast}$ contains all conjugates of
spaces of $K$.

The following theorem lists properties of these operations. In iteration of
the operations parentheses may be omitted.

Thus, $K^{\ast\ast}\overset{\operatorname*{def}}{=}\left(  (K)^{\ast}\right)
^{\ast}$; $HH(K)\overset{\operatorname*{def}}{=}H(H(K))$ and so on.

\begin{theorem}
Any set $K$ of finite dimensional Banach spaces has the following properties:

\begin{enumerate}
\item $K^{\ast\ast}=K$; $HH(K)=H(K)$; $QQ(K)=Q(K)$;

\item $K\subset H(K)$; $K\subset Q(K)$;

\item  If $K_{1}\subset K_{2}$ then $H(K_{1})\subset H(K_{2})$ and
$Q(K_{1})\subset Q(K_{2})$;

\item $\left(  H(K)\right)  ^{\ast}=Q(K^{\ast})$; $\left(  Q(K)\right)
^{\ast}=H(K^{\ast})$;

\item $HQ(HQ(K))=HQ(K)$; $QH(QH(K))=QH(K)$.
\end{enumerate}
\end{theorem}

\begin{proof}
1, 2 and 3 are obvious.

4. If $A\in Q(K)$ then $A=B/E$ for some $B\in K$ and its subspace $E$. So,
$A^{\ast}$ is isometric to a subspace of $B^{\ast}.$ Hence, $A^{\ast}\in
H(B^{\ast})$, i.e., $A^{\ast}\in H(K^{\ast})$. Since $A$ is arbitrary,
$\left(  Q(K)\right)  ^{\ast}\subseteq H(K^{\ast})$. Analogously, if $B\in K$
and $A\in H(B)$ then $A^{\ast}$ may be identified with a quotient $B^{\ast
}/A^{\perp}$, where $A^{\perp}$ is the annihilator of $A$ in $B^{\ast}$:%
\[
A^{\perp}=\{f\in B^{\ast}:f\left(  a\right)  =0\text{ \ for all }a\in A\}.
\]

Hence $A^{\ast}\in(Q(K^{\ast}))^{\ast}$ and thus $\left(  H(K)\right)  ^{\ast
}\subseteq Q(K^{\ast})$.

From the other hand,%
\begin{align*}
H(K^{\ast}) &  =\left(  H(K^{\ast})\right)  ^{\ast\ast}\subseteq\left(
Q(K^{\ast\ast})\right)  ^{\ast}=\left(  Q(K)\right)  ^{\ast};\\
Q(K^{\ast}) &  =\left(  Q(K^{\ast})\right)  ^{\ast\ast}\subseteq\left(
H(K^{\ast\ast})\right)  ^{\ast}=\left(  H(K)\right)  ^{\ast}.
\end{align*}

5. Let $A\in HQ(K)$. Then $A$ is isometric to a subspace of some quotient
space $E/F$, where $E\in K$; $F\hookrightarrow E$. If $B$ is a subspace of $A$
then $\left(  A/B\right)  ^{\ast}=\left(  E/F\right)  ^{\ast}/B^{\perp}$, i.e.
$\left(  Q\left(  HQ(K)\right)  \right)  ^{\ast}\subseteq Q\left(  Q(K)^{\ast
}\right)  $. Because of
\begin{align*}
Q\left(  HQ(K)\right)   &  =\left(  Q\left(  HQ(K)\right)  \right)  ^{\ast
\ast}\subseteq(Q\left(  Q(K)^{\ast}\right)  )^{\ast}\\
&  \subseteq H(Q(K)^{\ast\ast})=HQ(K),
\end{align*}

we have%
\[
H(Q(H(Q(K))))\subseteq H(H(Q(K)))=HQ(K).
\]

Analogously, if $A\in QH(K)$, then $A$ is isometric to a quotient space $F/E$,
where $F\in H(B)$ for some $B\in K$ and $E\hookrightarrow B$. If $W\in H(A)$,
i.e., if $W\in H(F/E)$ then $W^{\ast}=(F/E)^{\ast}/W^{\perp}$ and
$(F/E)^{\ast}$ is isometric to a subspace $E^{\perp}$ of $F^{\ast}%
\in(H(B))^{\ast}$. Thus, $(H(QH(K)))^{\ast}\subseteq H((H(K))^{\ast})$ and
\begin{align*}
H\left(  QH(K)\right)   &  =\left(  H\left(  QH(K)\right)  \right)  ^{\ast
\ast}\subseteq(H\left(  H(K)^{\ast}\right)  )^{\ast}\\
&  \subseteq Q(H(K)^{\ast\ast})=QH(K).
\end{align*}

Hence,
\[
Q(H(Q(HQ(K))))\subseteq Q(Q(H(K)))=QH(K).
\]

Converse inclusion follows from 2.
\end{proof}

It is obvious that for any class $W^{f}$ the set $\frak{N}=\frak{M}(W^{f})$
has following properties:

(\textbf{C}) $\frak{N}\ $\textit{is a closed subset of the Minkowski's space
}$\frak{M}$;

(\textbf{H}) \textit{If }$A\in\frak{N}$\textit{ and }$B\in H(A)$\textit{ then
}$B\in\frak{N}$;

(\textbf{A}$_{0}$)\textit{ For any }$A$\textit{, }$B\in\frak{N}$\textit{
\ there exists }$C\in\frak{N}$\textit{ \ such that }$A\in H(C)$\textit{ and}
$B\in H(C)$.

\begin{theorem}
Let $\frak{N}$ be a set of finite dimensional Banach spaces; $\frak{N}%
\subset\frak{M}$. If $\frak{N}$ has properties (\textbf{C}), (\textbf{H}) and
(\textbf{A}$_{0}$) then there exists a class $X^{f}$ such that $\frak{N}%
=\frak{M}(X^{f})$.
\end{theorem}

\begin{proof}
Conditions (\textbf{H}) and (\textbf{A}$_{0}$) in a natural (not unique) way
defines a partial order on $\frak{M}(X^{f})=\frak{N}$. If $D$ is an
ultrafilter which is consistent with this order, then the ultraproduct
$(\frak{N})_{D}$ of all spaces from $\frak{N}$ has desired properties.Since
the set $H\left(  (\frak{N})_{D}\right)  $ is closed, $\frak{M}(W^{f}%
)=\frak{N}$.
\end{proof}

Let $Y$ be a $B$-convex Banach space.

Let $X=l_{2}(Y)$ (and, hence, $X^{f}=\mathsf{D}_{2}(Y^{f})$). Consider the
Minkowski's base $\frak{M}(X^{f})$ and its enlargement $H(Q(\frak{M}%
(X^{f})))=HQ\frak{M}(X^{f})$.

\begin{theorem}
There exists a Banach space $W$ such that $HQ\frak{M}(X^{f})=\frak{M}(W^{f})$.
\end{theorem}

\begin{proof}
Obviously, $HQ\frak{M}(X^{f})$ has properties (\textbf{H}) and (\textbf{C}).

Since $\frak{M}(X^{f})$ is $2$-divisible, then for any $A,B\in\frak{M}(X^{f})$
the space $A\oplus_{2}B$ belongs to $\frak{M}(X^{f})$ and, hence, to
$HQ\frak{M}(X^{f})$. If $A,B\in Q\frak{M}(X^{f})$ then $A=F/F_{1}$;
$B=E/E_{1}$ for some $E,F\in Q\frak{M}(X^{f})$.

$F/F_{1}\oplus_{2}E$ is isometric to a space ($F\oplus_{2}E)/F_{1}^{\prime}$,
where
\[
F_{1}^{\prime}=\{(f,0)\in F\oplus_{2}E:f\in F_{1}\}
\]
and, hence, belongs to $Q\frak{M}(X^{f})$. Thus,
\[
F/F_{1}\oplus_{2}E/E_{1}=(F/F_{1}\oplus_{2}E)/E_{1}^{\prime},
\]
where
\[
E_{1}^{\prime}=\{(o,e)\in F/F_{1}\oplus_{2}E:e\in E_{1}\}
\]
and, hence, belongs to $Q\frak{M}(X^{f})$ as well.

If $A,B\in HQ\frak{M}(X^{f})$ then $A\hookrightarrow E$, $B\hookrightarrow F$
for some $E,F\in Q\frak{M}(X^{f})$.

$E\oplus_{2}F\in Q\frak{M}(X^{f})$ and, hence, $A\oplus_{2}B\in HQ\frak{M}%
(X^{f})$.

Thus, $HQ\frak{M}(X^{f})$ has the property (\textbf{A}$_{0}$). The desired
result follows from the preceding theorem.
\end{proof}

\begin{definition}
Let $X$ be a Banach space, which generates a class $X^{f}$\ of finite
equivalence. A class $\ast\ast(X^{f})$ is defined to be a class $W^{f}$ with
the Minkowski's base $\frak{M}(W^{f})=HQ\frak{M}(\mathsf{D}_{2}Y^{f})$.
\end{definition}

Clearly, $W^{f}$ is quotient closed. Obviously, $X^{f}<_{f}W^{f}$.

Let $\star$ be one more procedure that will be given by following steps.

Let $X\in\mathcal{B}$; $Y^{f}=\mathsf{D}_{2}\left(  X^{f}\right)  $. Let
$\left(  Y_{n}\right)  _{n<\infty}$ be a countable dense subset of
$\frak{M}(Y^{f})$. Consider the space $Z=(\sum_{n<\infty}\oplus Y_{n})_{2}$
and its conjugate $Z^{\ast}$.

$Z^{\ast}$ generates a class $\left(  Z^{\ast}\right)  ^{f}$, which will be
regarded as a result of the procedure $\star:X^{f}\rightarrow\left(  Z^{\ast
}\right)  ^{f}$. Iterations of the procedure $\star$ are given by following steps.

Let $\left(  Z_{n}\right)  _{n<\infty}$ be a countable dense subset of
$\frak{M}(\left(  Z^{\ast}\right)  ^{f})$. Consider a space $W=(\sum
_{n<\infty}\oplus Z_{n})_{2}$ and its conjugate $W^{\ast}$. Clearly,\ $W^{\ast
}$ generates a class $\left(  W^{\ast}\right)  ^{f}$, which may be regarded as
the result of the double procedure $\star$:
\[
\left(  W^{\ast}\right)  ^{f}=\star\left(  Z^{\ast}\right)  ^{f}=\star
\star\left(  X^{f}\right)  .
\]

\begin{theorem}
For any Banach space $X$ classes $\ast\ast\left(  X^{f}\right)  $ and
$\star\star\left(  X^{f}\right)  $ are identical.
\end{theorem}

\begin{proof}
From the construction it follows that $H(Z^{\ast})=\left(  QH(l_{2}(X)\right)
^{\ast}$ and that
\begin{align*}
H(W^{\ast}) &  =(QH(Z^{\ast}))^{\ast}=(Q(QH(l_{2}(X)))^{\ast})^{\ast}\\
&  =H(QH(l_{2}(X)))^{\ast\ast}=HQH(l_{2}(X)).
\end{align*}

Hence,
\[
\frak{M}(\left(  W^{\ast}\right)  ^{f})=HQ\frak{M}(Y^{f})=HQ\frak{M}%
(\mathsf{D}_{2}(X^{f})).
\]
\end{proof}

\begin{theorem}
Let $X$ be a Banach space, which generates a class of finite equivalence
$X^{f}$. If $X$ is of type $p>1$ and of cotype $q<\infty$, then the procedure
$\ast\ast$ maps $X^{f}$ to a class of the same type and cotype. If $X$ is
superreflexive then $\ast\ast\left(  X^{f}\right)  $ is superreflexive too.
\end{theorem}

\begin{proof}
According to [9], for any $B$-convex Banach space $X$ of nontrivial type $p$
(and, hence, cotype $q$) its conjugate $X^{\ast}$ is of type $p^{\prime
}=q/(q-1)$ and of cotype $q^{\prime}=p/(p-1)$ (and, hence, is $B$ convex as
well). Obviously, $\star\left(  X^{f}\right)  $ is of type $p^{\prime
}=q/(q-1)$ and of cotype $q^{\prime}=p/(p-1)$ and thus $\star\star\left(
X^{f}\right)  =\ast\ast\left(  X^{f}\right)  $ is of type $p$ and of
cotype\textit{ }$q$. It is clear that the procedure $\star$ holds
superreflexivity. From the preceding theorem it follows that $\ast\ast$ also
has this property.
\end{proof}

\begin{remark}
If $X$ is not $B$-convex, then
\[
\star\left(  X^{f}\right)  =\star\star\left(  X^{f}\right)  =\ast\ast\left(
X^{f}\right)  =\left(  c_{0}\right)  ^{f}.
\]
\end{remark}

\begin{theorem}
For any Banach space $X$ the class $\ast\ast\left(  X^{f}\right)  $ is $2$-divisible.
\end{theorem}

\begin{proof}
Let $\frak{N}=\frak{M}(\ast\ast(X^{f}))$.

Since for any pair $A,B\in\frak{N}$ \ their $l_{2}$-sum belongs to $\frak{N}$,
then, by induction, $(\sum\nolimits_{i\in I}\oplus A_{i})_{2}\in\frak{N}$\ for
any finite subset $\{A_{i}:i\in I\}\subset\frak{N}$.

Hence any infinite direct $l_{2}$-sum $(\sum\nolimits_{i\in I}\oplus
A_{i})_{2}$ is finite representable in $\ast\ast\left(  X^{f}\right)  $.

Let $\{A_{i}:i<\infty\}\subset\frak{N}$ \ be dense in $\frak{N}$.

Let $Y_{1}=(\sum\nolimits_{i<\infty}\oplus A_{i})_{2}$; $Y_{n+1}=Y_{n}%
\oplus_{2}Y_{1}$; $Y_{\infty}=\overline{\cup Y_{n}}$ (the upper line marks the closure).

Clearly, $Y_{\infty}=l_{2}\left(  Y_{\infty}\right)  $ belongs to $\ast
\ast\left(  X^{f}\right)  $.
\end{proof}

\section{Almost $\omega$-homogeneous Banach spaces}

\begin{definition}
Let $X$ be a Banach space; $\mathcal{K}$ be a class of Banach spaces. $X$ is
said to be almost $\omega$-homogeneous with respect to $\mathcal{K}$ if for
any pair of spaces $A$, $B$ of $\mathcal{K}$ such that $A$ is a subspace of
$B$ ($A\hookrightarrow B$), every $\varepsilon>0$ and every isometric
embedding $i:A\rightarrow X$ there exists an isomorphic embedding
$\hat{\imath}:B\rightarrow X$, which extends $i$ (i.e., $\hat{\imath}|_{A}=i$)
such that
\[
\left\|  \hat{\imath}\right\|  \left\|  \hat{\imath}^{-1}\right\|
\leq(1+\varepsilon).
\]
\end{definition}

If $X$ is almost $\omega$-homogeneous with respect to $H(X)$ it will be
referred to as an almost $\omega$-homogeneous space.

The main role in the next plays the notion of the\textit{\ amalgamation}
property. It will be convenient to introduce some terminology.

A fifth $v=\left\langle A,B_{1},B_{2},i_{1},i_{2}\right\rangle $, where $A$,
$B_{1}$, $B_{2}\in\frak{M}(X^{f})$; $i_{1}:A\rightarrow B_{1}$ and$\ i_{2}%
:A\rightarrow B_{2}$ are isometric embeddings, will be called \textit{the }%
$V$\textit{-formation over} $\frak{M}(X^{f})$. The space $A$ will be called
\textit{the root of the }$V$\textit{-formation }$v$\textit{. }If there
exists\textit{\ }a triple $t=\left\langle j_{1},j_{2},F\right\rangle $ so that
$F\in\frak{M}(X^{f})$; $\ j_{1}:B_{1}\rightarrow F$ and $j_{2}:B_{2}%
\rightarrow F$ are isometric embeddings such that $j_{1}\circ i_{1}=j_{2}\circ
i_{2}$, then the $V$-formation $v$ is said to be \textit{amalgamated in}
$\frak{M}(X^{f})$, and the triple $t$ (or simply the space $F$) is said to be
its \textit{amalgam}.

Let $\operatorname{Amalg}(\frak{M}(X^{f}))$ be a set of all spaces
$A\in\frak{M}(X^{f})$ with the property:

\textit{Any }$V$\textit{-formation }$v$\textit{ with root }$A$\textit{\ is
amalgamated in }$\frak{M}(X^{f})$.

\begin{definition}
Let $X\in\mathcal{B}$ \ generates a class $X^{f}$\ with the Minkowski's base
$\frak{M}(X^{f})$. It will be said that $\frak{M}(X^{f})$ (and the whole class
$X^{f}$) has the amalgamation property if
\[
\frak{M}(X^{f})=\operatorname{Amalg}(\frak{M}(X^{f}))
\]
\end{definition}

\begin{theorem}
For any class $X^{f}$ having the amalgamation property, $E\in\mathcal{E}%
\left(  X^{f}\right)  $ if and only if $E\in X^{f}$ and $E$ is almost $\omega$-homogeneous.
\end{theorem}

\begin{proof}
Let $E\hookrightarrow Z\in X^{f}$ and $E$ be almost $\omega$-homogeneous. Let
$A\hookrightarrow Z$ be a finite dimensional subspace; $E\cap A=B$ and
$\varepsilon>0$. Consider the identical embedding $id_{B}:B\rightarrow E$.
Since $B\hookrightarrow A$, $id_{B}$ may be extended to an embedding
$u:A\rightarrow E$ with $\left\|  u\right\|  \left\|  u^{-1}\right\|
\leq1+\varepsilon$. Thus, $E\in\mathcal{E}\left(  X^{f}\right)  $.

Conversely, let $E\in\mathcal{E}\left(  X^{f}\right)  $; $A\hookrightarrow E$
be finite dimensional. Let $i:A\rightarrow E$ be an operator. Let
$iA\hookrightarrow B\in\frak{M}\left(  X^{f}\right)  $ for some $B$.

Consider a space $Z$ such that $E\hookrightarrow Z$; $iA\hookrightarrow
B\hookrightarrow Z$. Such space exists because of the amalgamation property of
$\frak{M}\left(  X^{f}\right)  $. Since $E\in\mathcal{E}\left(  X^{f}\right)
$, $E\prec_{u}Z$, i.e. there is an embedding $u:B\rightarrow E$ such that
$\left\|  u\right\|  \left\|  u^{-1}\right\|  \leq1+\varepsilon$, which is
identical on the intersection $E\cap B$; $u\mid_{E\cap B}=Id_{E\cap B}$. Since
$iA\hookrightarrow B$ and $iA\hookrightarrow E$, $iA\hookrightarrow E\cap B$.
Clearly, $u$ extends the embedding $i:A\rightarrow B$. Since
$iA\hookrightarrow B$ and $\varepsilon$ are arbitrary, $E$ is almost $\omega
$-homogeneous .
\end{proof}

The one-side shuttle procedure shows that $E$ is almost universal for all
separable spaces that are finitely representable in $X^{f}$ (in other
terminology, any separable $E\in\mathcal{E}\left(  X^{f}\right)  $ is an
\textit{approximative envelope} of the class $X^{f}$).

\begin{theorem}
Let $X^{f}$ has the amalgamation property; $E\in\mathcal{E}\left(
X^{f}\right)  $ be separable. For any separable space $Z<_{f}E$ and every
$\varepsilon>0$ there exists an isomorphic embedding $u:Z\rightarrow E$ with
$\left\|  u\right\|  \left\|  u^{-1}\right\|  \leq1+\varepsilon$.
\end{theorem}

\begin{proof}
Let $Z=\overline{\cup Z_{n}}$, where $Z_{1}\hookrightarrow Z_{2}%
\hookrightarrow...$ be an increasing chain of finite dimensional subspaces,
starting with 1-dimensional space $Z_{1}$. Define inductively a sequence of
isomorphic embeddings $\left(  i_{n}\right)  $. Let $i_{1}:Z_{1}\rightarrow
E_{X}$ be an (isometric) embedding. Let $i_{n+1}:Z_{n+1}\rightarrow Z$ be an
$\left(  1+\varepsilon^{n}\right)  $-isomorphic embedding, which extends
$i_{n}$. Certainly, $\overline{\cup i_{n}Z_{n}}=Z^{\prime}\hookrightarrow
E_{X}$ is $\left(  1+2\varepsilon\right)  $-isomorphic to $Z$.
\end{proof}

\begin{remark}
The same result may be obtained in other way. It may be shown that\ a
separable almost $\omega$-homogeneous space $E$ (which is existentialy closed)
is unique up to almost isometry.

Since, according to $[8]$, \textbf{every }separable space\textbf{ }$Z<_{f}E$
may be isometricaly embedded into \textbf{some }separable\textbf{ }%
$E\in\mathcal{E}\left(  E^{f}\right)  $, immediately $E$ has the desired property.
\end{remark}

Now it will be shown that every class $X^{f}$ of finite equivalence, which is
divisible and quotient-closed is crudely finite equivalent to a class $\Gamma
X^{f}$, which enjoys the isomorphic amalgamation property.

To define the procedure $\Gamma$, which sends a divisible quotient-closed
class $X^{f}$ to the desired class $\Gamma X^{f}$, which will be called
\textit{the Gurarii compression of }$X^{f}$ consider a space $Y\in X^{f}$ with
$H(Y)=\frak{M}(X^{f})$ (recall that every ultrapower $\left(  Z\right)  _{D}$
of any $Z\in X^{f}$ has the such property). By the theorem 7, $Y\oplus Y$
belongs to $X^{f}$ for an orthogonal direct sum $\oplus$. So, for any pair
$A$, $B\in\frak{M}(X^{f})$ is defined their orthogonal direct sum $A\oplus B$
in a determined way (chose $A$, resp. $B$ as a subspace of the first, resp.,
of the second component of the sum $Y\oplus Y$). Let $\left\|  \left(
a,b\right)  \right\|  $ be the corresponding norm on $Y\oplus Y$ (and, hence,
on $A\oplus B$).

Consider $\frak{M}(X^{f})$ and for every $A\in\frak{M}(X^{f})$ define a new
norm, say, $\left|  \left\|  \cdot\right\|  \right|  $, by the rule: for $a\in
A$%
\[
\left|  \left\|  a\right\|  \right|  =\inf\{\left\|  \left(  u,v\right)
\right\|  :u+v=a\}
\]
Let us show that this is really the norm.

\begin{proposition}
The function $\left|  \left\|  \cdot\right\|  \right|  :a\rightarrow\left|
\left\|  a\right\|  \right|  $ is positive, homogeneous and satisfies the
triangle inequality.
\end{proposition}

\begin{proof}
Immediately, $\left|  \left\|  a\right\|  \right|  \geq0$: $\left|  \left\|
a\right\|  \right|  =0\Leftrightarrow a=0$ and, for a scalar $\lambda$,
$\left|  \left\|  \lambda a\right\|  \right|  =$ $\left|  \lambda\right|
\left|  \left\|  a\right\|  \right|  $. Prove the triangle inequality.
\begin{align*}
\left|  \left\|  a+b\right\|  \right|   &  =\inf\{\left\|  \left(  u,v\right)
\right\|  :u+v=a+b\}\\
&  =\inf\{\left\|  \left(  u_{1}+v_{1},u_{2}+v_{2}\right)  \right\|
:u_{1}+v_{1}+u_{2}+v_{2}=a+b\}\\
&  =\inf\{\left\|  \left(  u_{1},u_{2}\right)  +\left(  v_{1}+v_{2}\right)
\right\|  :u_{1}+v_{1}+u_{2}+v_{2}=a+b\}\\
&  \leq\inf\{\left\|  \left(  u_{1},u_{2}\right)  \right\|  +\left\|  \left(
v_{1}+v_{2}\right)  \right\|  :u_{1}+v_{1}+u_{2}+v_{2}=a+b\}\\
&  \leq\inf\{\left\|  \left(  u_{1},u_{2}\right)  \right\|  :u_{1}%
+u_{2}=a\}+\inf\{\left\|  \left(  v_{1},v_{2}\right)  \right\|  :v_{1}%
+v_{2}=b\}\\
&  \leq\left|  \left\|  a\right\|  \right|  +\left|  \left\|  b\right\|
\right|  .
\end{align*}
\end{proof}

\begin{theorem}
Norms $\left|  \left\|  a\right\|  \right|  $ and $\left\|  a\right\|  $ are equivalent.
\end{theorem}

\begin{proof}
Immediately, $\left|  \left\|  a\right\|  \right|  \leq\left\|  a\right\|  $
for all $a\in A$. The inverse inequality easily follows from another
definition of the norm $\left|  \left\|  \cdot\right\|  \right|  $, which may
be identified with the norm of the quotient $\left(  A\oplus A\right)  /N$,
where
\[
N=\{\left(  v,-v\right)  :v\in A\}\hookrightarrow A\oplus A.
\]

Indeed, for a class $\left[  \left(  a,0\right)  \right]  $, generated by an
element $\left(  x,0\right)  \in A\oplus A$ its quotient norm is
\begin{align*}
\left\|  \left(  a,0\right)  \right\|  ^{Q} &  =\inf\{\left\|  \left(
a+v,-v\right)  \right\|  :v\in A\}\geq\inf\{\max\{\left\|  a+v\right\|
,\left\|  v\right\|  \}:v\in A\}\\
&  \geq\inf\{\max\{\left|  \left\|  a\right\|  -\left\|  v\right\|  \right|
,\left\|  v\right\|  \}:v\in A\}\\
&  =\left\|  a\right\|  \inf\{\max\{\left|  1-\lambda\right|  ,\left|
\lambda\right|  \}:\lambda\in\mathbb{R}\}\geq1/2\left\|  a\right\|
\end{align*}
and is equal to%
\begin{align*}
\left\|  \left(  a,0\right)  \right\|  ^{Q} &  =\inf\{\left\|  \left(
a+v,-v\right)  \right\|  :v\in A\}\\
&  =\inf\{\left\|  \left(  u_{1},u_{2}\right)  :u_{1}+u_{2}=a\right\|
=\left|  \left\|  a\right\|  \right|  .
\end{align*}

So,
\[
1/2\left\|  a\right\|  \leq\left|  \left\|  a\right\|  \right|  \leq\left\|
a\right\|  .
\]
\end{proof}

Now define the procedure $\Gamma$ of Gurarii compression.

Let $\mathcal{K}$ be a class of finite-dimensional Banach spaces. For every
$A\in\mathcal{K}$ put $R(A)=\left\langle A,\left|  \left\|  \cdot\right\|
\right|  \right\rangle $, where $\left|  \left\|  \cdot\right\|  \right|  $ is
defined as above. Put%
\[
\Gamma\left(  \mathcal{K}\right)  =\{R(A):A\in\mathcal{K}\}.
\]

Let $X^{f}$ be a quotient-closed divisible class; $\frak{M}(X^{f})$ be its
Minkowski's base. Consider a set $\Gamma(\frak{M}(X^{f})).$ Since $R(A)$ may
be identified with the quotient $\left(  A\oplus A\right)  /N$ where
$N\in\frak{M}(X^{f})$ and because of $\frak{M}(X^{f})$ is quotient-closed, it
follows that%
\[
\Gamma(\frak{M}(X^{f}))\subseteq\frak{M}(X^{f}).
\]

It is obvious that $R(A\oplus B)=R(A)\oplus R(B).$ This shows that the set
$\Gamma(\frak{M}(X^{f}))$ enjoys the property (\textbf{A}$_{0}).$ Properties
(\textbf{H}) and (\textbf{C}) for $\Gamma(\frak{M}(X^{f}))$ are obvious as
well. So, $\Gamma(\frak{M}(X^{f}))$ may be regarded as a Minkowski's base of a
class $W^{f}\overset{\operatorname{def}}{=}\Gamma(X^{f})$. It is clear that
this procedure is idempotent, i.e. $\Gamma\Gamma(X^{f})=\Gamma(X^{f})$.
Besides, $\Gamma(X^{f})$ is quotient-closed. The class $\Gamma(X^{f})$ will be
referred to as \textit{the Gurarii compression of the class }$X^{f}$.

\begin{theorem}
Let $X^{f}$ be quotient-closed divisible class. Its \textit{Gurarii
compression }$\Gamma(X^{f})$ is finitely representable in $X^{f}$ and is
crudely finitely equivalent to $X^{f}$. The class $\Gamma(X^{f})$ is
quotient-closed divisible as well.
\end{theorem}

\begin{proof}
The first part follows from the inclusion $\Gamma(X^{f})\subseteq
\frak{M}(X^{f})$. The second one is a consequence of the previous theorem. The
third pat is obvious.
\end{proof}

\begin{theorem}
Let $X^{f}$ be a quotient closed, \textit{divisible class. }Its
\textit{Gurarii compression }$W^{f}=\Gamma(X^{f})$ enjoys the  amalgamation property.
\end{theorem}

\begin{proof}
Let $A$, $B_{1}$, $B_{2}\in\frak{M}(W^{f})$; $i_{1}:A\rightarrow B_{1}$ and
$i_{2}:A\rightarrow B_{2}$ be isometric embeddings.

Assume for simplicity that $A\hookrightarrow B_{1}$ and $A\hookrightarrow
B_{2}$.

Consider a linear space $F$ formed by elements $f$ of kind $f=x+y$, where
$x\in B_{1}$ and $y\in B_{2}$, which is equipped with the norm:
\[
\left\|  f\right\|  _{F}=\inf\{\sqrt{\left\|  x\right\|  _{B_{1}}^{2}+\left\|
y\right\|  _{B_{2}}^{2}}:x+y=f\}.
\]

Obviously, the space $F$ is an amalgam of the $V$-formation $\ \left\langle
A,B_{1},B_{2},i_{1},i_{2}\right\rangle $.

To close the proof notice that the space $F$ may be identified with the
quotient $E/H$, where $E=B_{1}\oplus_{2}B_{2}\ $and\ $H$ is a subspace of $E$,%
\[
H=\{\left(  a,-a\right)  :a\in A\}.
\]

Surely, $E/H\in\frak{M}(W^{f})$ and, hence, $F\in\frak{M}(W^{f})$ as well.
\end{proof}

\section{A solution of Bourgain's problem.}

\begin{theorem}
There exists a separable reflexive Banach space $G$ with the following
property: for any separable superreflexive Banach space $X$ and each
$\varepsilon>0$ there exists an isomorphic embedding $u:X\rightarrow G$ such
that $\left\|  u\right\|  \left\|  u^{-1}\right\|  \leq1+\varepsilon$.
\end{theorem}

\begin{proof}
Let $1<p\leq2\leq q<\infty$; $C(p,q,N)$ $\subset\mathcal{B}$ be a set of all
such separable superreflexive Banach spaces $X$ that are both of type $p$ and
of cotype $q$ with $\max\{T_{p}\left(  X\right)  ,C_{q}\left(  X\right)
\}\leq N$.

Let $n\in\mathbb{N}$;$\ p=1-1/n$; $q=n$; $N=n$, $C_{n}=C(1-1/n,n,n)$.

Consider a Banach space $X_{n}=(\sum\oplus\{E:E\in C_{n}\})_{2}$ and a class
$\left(  X_{n}\right)  ^{f}$ that is generated by $X_{n}$.

The procedure $\ast\ast$ sends the class $\left(  X_{n}\right)  ^{f}$ to the
class $\left(  W_{n}\right)  ^{f}=\ast\ast((X_{n})^{f})$, which is
superreflexive. Its Gurarii compression - the class $\Gamma\left(
W_{n}\right)  ^{f}$ has the amalgamation property and hence contains a
separable approximative envelope, say, $E_{n}$. Surely, each space from
$C(1-1/n,n,n)$ is isomorphic to some subspace of $E_{n}$.

Consider a space
\[
G=(\sum\nolimits_{n=1}^{\infty}\oplus E_{n})_{2}%
\]

Since a class $\mathcal{SSR}$ of all separable superreflexive Banach spaces
may be represented as the union
\[
\mathcal{SSR}=\cup\{C_{n}:n<\infty\},
\]
it is clear that the space $G$ is desired.
\end{proof}

\begin{remark}
A space $G$ with the property of the preceding theorem is not unique. Any
space that is an $l_{p}$-sum of $E_{n}$'s has the same property. Moreover, if
$F\in\left(  W_{n}\right)  ^{f}$ is separable then $F\oplus_{2}E_{n}$ is also
an approximative separable envelope of $\left(  W_{n}\right)  ^{f}$. It is
clear that there exists a continuum number of pairwice non isomorphic spaces
of kind $F\oplus_{2}E_{n}$.
\end{remark}

Recall that a Banach space $X$ is said to be \textit{complementably universal}
for a class $\mathcal{K}$ of Banach spaces if every $E\in\mathcal{K}$ may be
isomorphicaly embedded into $X$ in a such way that its isomorphic copy
$E^{\prime}$ in $X$ admits a (bounded linear) projection $P:X\rightarrow E$.

From $[10]$ it follows that there no separable (non necessary reflexive)
Banach space is complementary universal for the class of all subspaces of
$l_{p}$ for any given $p$, $1\leq p\neq2<\infty$. Hence there not exists a
complementably universal separable space for the whole class $\mathcal{SSR}$.
At the same time, according to [2], there is a reflexive Banach space which is
complementably universal for a class of those spaces from $\mathcal{SSR}$ that
have the (metric) approximation property.

\begin{acknowledgement}
Author wish to express his gratitude to Hanfeng Li, who pointed out on a gap
in the first version of this paper for valuable discussions and interest to
this work.
\end{acknowledgement}

\section{References}

\begin{enumerate}
\item  Bourgain J. \textit{On separable Banach spaces universal for all
separable reflexive spaces}, Proc. AMS \textbf{79:2} (1980) 241-246

\item  Prus S. \textit{Finite - dimensional decompositions with p - estimates
and universal Banach spaces}, Bull. Polish. Acad. Sci. Math.\textbf{ 31:5-8}
(1983) 281-288

\item  Heinrich S. \textit{The isomorphic problem of envelops}, Studia Math.
\textbf{73:1 }(1982) 41-49

\item  Schwartz L. \textit{Geometry and probability in Banach spaces}, Bull.
AMS \textbf{4:2} (1981) 135-141

\item  Maurey B. and Pisier G. \textit{S\'{e}ries de variables al\'{e}atoires
vectorielles ind\'{e}pen-dantes et propri\'{e}t\'{e}s g\'{e}om\'{e}triques des
espaces de Banach}, Studia Math. \textbf{58 }(1976) 45-90

\item  Tokarev E.V. \textit{Spectrum of infinite-dimensional Banach spaces
}(transl. from Russian), Functional Analysis and its Applications \textbf{23},
No.1 (1989) 76-78

\item  Stern J. \textit{Ultrapowers and local properties in Banach spaces},
Trans. AMS \textbf{240} (1978) 231-252

\item  Tokarev E.V. \textit{Injective Banach spaces in the finite equivalence
classes }(transl. from Russian), Ukrainian Mathematical Journal \textbf{39:6}
(1987) 614-619

\item  Pisier G. \textit{On the duality between type and cotype}, Lect. Notes
in Math. \textbf{939} (1982)

\item  Johnson W.B., Szankowski A. \textit{Complementably universal Banach
spaces}, Studia Math. \textbf{58} (1976) 91-97
\end{enumerate}
\end{document}